\documentclass[a4paper,12pt]{amsart}
\usepackage{pslatex}
\usepackage{mathrsfs}
\usepackage{dsfont}
\usepackage[colorlinks]{hyperref}
\usepackage[a4paper,asymmetric]{geometry}
\usepackage{pgf}
\newtheorem{theorem}{Theorem}[section]

\newtheorem{lemma}[theorem]{Lemma}
\newtheorem{proposition}[theorem]{Proposition}
\theoremstyle{definition}
\newtheorem{definition}[theorem]{Definition}
\theoremstyle{remark}
\newtheorem{remark}[theorem]{Remark}

\numberwithin{equation}{section}
\newcommand{\arxiv}[2]{\href{http://arxiv.org/abs/#2}{arxiv:#1/#2}}
\DeclareMathOperator{\Div}{div}
\DeclareMathOperator{\Tr}{Tr}
\DeclareMathOperator{\e}{e}
\newcommand{\Erre}{\mathbf{R}}
\newcommand{\N}{\mathbf{N}}
\newcommand{\Pb}{\mathbb{P}}
\newcommand{\E}{\mathbb{E}}
\newcommand{\W}{\mathcal{W}}
\newcommand{\ris}{\mathcal{R}}
\newcommand{\semi}{\mathcal{P}}
\newcommand{\gen}{\mathscr{L}}
\newcommand{\genG}{\mathscr{L}^{\textsf{\Tiny G}}}
\newcommand{\genSV}{\mathscr{L}_{\textsf{\Tiny SY}}}
\newcommand{\ms}{\mathscr{C}}
\newcommand{\msP}{\mathscr{C}_{\textsf{\Tiny P}}}

\newcommand{\msSV}{\mathscr{C}_{\textsf{\Tiny SY}}}
\newcommand{\B}{\mathscr{B}}
\newcommand{\F}{\mathscr{F}}
\newcommand{\bm}{\mathcal{B}_b}
\newcommand{\eps}{\epsilon}
\newcommand{\er}{\text{\Tiny$(R)$}}
\newcommand{\ern}{\text{\Tiny$(R_n)$}}
\newcommand{\Torus}{{\mathds{T}_\text{\tiny$3$}}}
\newcommand{\cov}{\mathcal{Q}}
\newcommand{\loc}{\mathrm{loc}}
\newcommand{\scal}[1]{\langle #1\rangle}
\newcommand{\uno}{\mathds{1}}
\newcommand{\Energy}{\mathscr{E}}
\newcommand{\tX}{\widetilde{X}}
\newcommand{\Xz}{\widehat{X}}
\newcommand{\curva}[1]{\begin{pgftranslate}{\pgfxy(#1)}\soluzione\pgfusepath{stroke}\end{pgftranslate}}
\def\soluzione{%
  \pgfpathmoveto{\pgfxy( 0,  0  )}  \pgfpathlineto{\pgfxy( 4,  1.2)}
  \pgfpathlineto{\pgfxy( 8,  3.9)}  \pgfpathlineto{\pgfxy(12,  7  )}
  \pgfpathlineto{\pgfxy(16, 10  )}  \pgfpathlineto{\pgfxy(20, 12.7)}
  \pgfpathlineto{\pgfxy(24, 15.1)}  \pgfpathlineto{\pgfxy(28, 17.1)}
  \pgfpathlineto{\pgfxy(32, 18.7)}  \pgfpathlineto{\pgfxy(36, 20  )}
  \pgfpathlineto{\pgfxy(40, 21.1)}  \pgfpathlineto{\pgfxy(44, 21.9)}
  \pgfpathlineto{\pgfxy(48, 22.6)}  \pgfpathlineto{\pgfxy(52, 23.1)}
  \pgfpathlineto{\pgfxy(56, 23.5)}  \pgfpathlineto{\pgfxy(60, 23.8)}
  \pgfpathlineto{\pgfxy(64, 24.1)}  \pgfpathlineto{\pgfxy(68, 24.3)}
  \pgfpathlineto{\pgfxy(72, 24.4)}  \pgfpathlineto{\pgfxy(76, 24.6)}
  \pgfpathlineto{\pgfxy(80, 24.7)}
}
\begin{document}
\title[The martingale problem for the Navier-Stokes equations]{The martingale problem for Markov solutions to the Navier-Stokes equations}
\author[M. Romito]{Marco Romito}
\address{Dipartimento di Matematica, Universit\`a di Firenze
viale Morgagni 67/A, I-50134, Firenze, Italia
\url{http://www.math.unifi.it/users/romito}}
\email{romito@math.unifi.it}
\date{}
\begin{abstract}
Under suitable assumptions of regularity and non-degeneracy on the covariance
of the driving additive noise, any Markov solution to the stochastic Navier-Stokes
equations has an associated generator of the diffusion and is the unique solution
to the corresponding martingale problem. Some elementary examples are discussed
to interpret these results.
\end{abstract}
\maketitle
\section{Introduction}

For some interesting stochastic partial differential equations, such as
the three dimensional Navier-Stokes equations, well-posedness of the
associated martingale problem is still an open problem. If on one hand
this corresponds to a poor understanding of the deterministic dynamics
(see for instance Feffermann~\cite{Fef00}), on the other hand there is
still the possibility that the stochastic problem may have better
uniqueness properties as in the finite dimensional case, see for example
\cite[Chapter 8]{StrVar79}, see also~\cite{Fla08} for a review in comparison
with the infinite dimension and~\cite{FlaGubPri08} for a infinite dimensional
positive result.

With these uniqueness problems in mind, it is then reasonable, and sometimes
useful, to consider special solutions with additional properties. We focus
here on solutions to the stochastic Navier-Stokes equations,
\begin{equation}\label{e:nse}
\begin{cases}
\dot u - \nu\Delta u + (u\cdot\nabla)u + \nabla p = \eta,\\
\Div u = 0,
\end{cases}
\end{equation}
which are Markov processes, and we give a short account of the theory
introduced in~\cite{FlaRom06} and \cite{FlaRom08}, and developed
in~\cite{FlaRom07}, \cite{Rom08} (see also~\cite{BloFlaRom08}, where
such ideas have been applied on a stochastic fourth order parabolic
equation driven by space-time white noise and modelling surface growth).
It turns out that, under suitable assumptions of regularity and
non-degeneracy on the covariance of the noise, all Markov solutions
to~\eqref{e:nse} are strong Feller and converge, exponentially fast,
to a unique invariant measure. Similar (and sometimes better) results
have been already obtained by Da Prato \& Debussche~\cite{DapDeb03}
(see also \cite{DebOda06}, \cite{Oda07} and \cite{DapDeb08}) with
a completely different method.

The aim of this paper is to understand the martingale problem associated
to any of the Markov processes which are solutions to~\eqref{e:nse}.
We prove that each of these Markov processes has a generator and it is
the unique solution of the martingale problem associated to the
generator of its own dynamics. It remains completely open to understand
what these generators have in common with the formal generator
$$
\gen^\star
 = \frac12\Tr[\cov D^2] - \scal{-\nu\nabla u + (u\cdot\nabla), D}.
$$
where $\cov$ is the spatial covariance of the noise, and in Section~\ref{ss:generator}
we make an attempt to connect $\gen^\star$ with the generators.

In order to compare all such results, in the final part of the paper
we give a few instructive examples (essentially taken from Stroock
\& Yor~\cite{StrYor80}) of one dimensional stochastic differential
equations where uniqueness is not ensured. All these examples show
that it is possible to have problems where different solutions are
Feller (and each of them has a unique invariant measure, see
Section~\ref{ss:girsanov}) or even strong Feller
(Section~\ref{ss:stroockyor}).

In conclusion, regardless of the improvement gained with the
Markov approach (\cite{DapDeb03} and \cite{FlaRom08}), the
problem remains completely open. The examples presented here
strengthen the belief that we must understand the behaviour
of these solutions when they approach \emph{blow-up} epochs
(see for example \cite{FlaRom02} for an attempt in this direction).

The paper is organised as follows. Section~\ref{s:markov} contains
a short summary of results on Markov solutions for the Navier-Stokes
equations. Existence of the generator and all details on the
martingale problem can be found in Section~\ref{s:mp}. Finally,
the examples are presented in the last section.
\section{Markov solutions for the Navier-Stokes equations}\label{s:markov}

In this section we give a short account of the ideas introduced in \cite{FlaRom06},
\cite{FlaRom07}, \cite{FlaRom08}, \cite{Rom08}. We focus on the equations on the
three-dimensional torus $\Torus=[0,2\pi]^3$ with periodic boundary conditions.

We fix some notations we shall use throughout the paper and we refer to
Temam \cite{Tem83} for a detailed account of all the definitions. Let
$\mathcal{D}^\infty$ be the space of infinitely differentiable
divergence-free periodic vector fields $\varphi:\Erre^3\to\Erre^3$
with mean zero on $\Torus$. Denote by $H$ the closure of $\mathcal{D}^\infty$
in $L^2(\Torus,\Erre^3)$ and by $V$ the closure in $H^1(\Torus,\Erre^3)$.
Denote by $A$, with domain $D(A)$, the \emph{Stokes} operator and
define the bi-linear operator $B:V\times V\to V'$ as the projection
onto $H$ of the nonlinearity of equation~\eqref{e:nse}.
Let $\cov$ be a linear bounded symmetric positive operator on $H$ with finite trace.
Consider finally the abstract form of problem~\ref{e:nse},
\begin{equation}\label{e:abstractnse}
du + \nu Au + B(u,u) = \cov^{\frac12}dW,
\end{equation}
where $W$ is a cylindrical Wiener process on $H$.

The first step is to define a solution to the above equation~\eqref{e:abstractnse}.
To this end, let $\Omega=C([0,\infty);D(A)')$, let $\B$ be the Borel $\sigma$-field
on $\Omega$ and let $\xi:\Omega\to D(A)'$ be the canonical process on
$\Omega$ (that is, $\xi_t(\omega)=\omega(t)$). A filtration can be
defined on $\B$ as $\B_t=\sigma(\xi_s:0\leq s\leq t)$.

For every $\varphi\in\mathcal{D}^\infty$ consider the process $(M_t^\varphi)_{t\geq0}$
on $\Omega$ defined for $t\geq0$ as
\begin{equation}\label{e:mart_eq}
M_t^\varphi
= \scal{\xi_t-\xi_0,\varphi}
 +\nu\int_0^t\scal{\xi_s,A\varphi}\,ds
 -\int_0^t\scal{B(\xi_s,\varphi),\xi_s}\,ds.
\end{equation}
Moreover, for every $n\in\N$, $n\geq1$, define for $t\geq0$ the process
\begin{equation}\label{e:mart_energy}
\Energy_t^n
= |\xi_t|^{2n}_H
 +2n\nu\int_0^t|\xi_s|^{2n-2}_H|\xi_s|_V^2\,ds
 -n(2n-1)\Tr[\cov]\int_0^t|\xi_s|_H^{2n-2}\,ds,
\end{equation}
when $\xi\in L_\loc^\infty([0,\infty);H)\cap L^2_\loc([0,\infty);V)$,
and $\infty$ elsewhere.
\begin{definition}\label{d:ems}
Given $\mu_0\in\Pr(H)$, a probability $P$ on $(\Omega,\B)$ is a solution
starting at $\mu_0$ to the martingale problem associated to the Navier-Stokes
equations \eqref{e:abstractnse} if
\begin{itemize}
\item $P[L_\loc^\infty([0,\infty);H)\cap L^2_\loc([0,\infty);V)]=1$,
\item for each $\varphi\in\mathcal{D}^\infty$ the process $M_t^\varphi$
is square integrable and $(M_t^\varphi,\B_t,P)$ is a continuous martingale
with quadratic variation $[M^\varphi]_t=t|\cov^{\frac12}\varphi|^2_H$,
\item for each $n\geq1$, the process $\Energy^n_t$ is $P$-integrable and
for almost every $s\geq 0$ (including $s=0$) and all $t\geq s$,
$$
\E[\Energy_t^n|\B_s]\leq \Energy_s^n,
$$
\item the marginal of $P$ at time $t=0$ is $\mu_0$.
\end{itemize}
\end{definition}
\begin{remark}
Goldys, Rockner \& Zhang~\cite{GolRocZha08} have pointed out that,
due to a lack of measurability of conditional probabilities, the
condition on the process $\Energy$ should be replaced by an analogous
condition on $\sup_{s\leq t}\Energy_s$.
\end{remark}
The idea behind the existence of Markov solutions is by Krylov~\cite{Kry73}
(see also Chapter 12 of Stroock \& Varadhan~\cite{StrVar79}). Define
for every $x\in H$ the set $\ms(x)$ of all solutions (according to
Definition~\ref{d:ems} above) starting at $\delta_x$.
\begin{theorem}[Theorem $4.1$, \cite{FlaRom08}]\label{t:markov}
There exists a family $(P_x)_{x\in H}$ of weak martingale solutions
such that $P_x\in\ms(x)$ for every $x\in H$ and the \emph{almost sure
Markov property} holds. More precisely, for every $x\in H$, for almost
every $s\geq0$ (including $s=0$), for all $t\geq s$ and all bounded
measurable $\phi:H\to\Erre$,
\begin{equation}\label{e:markov}
\E^{P_x}[\phi(\xi_t')|\B_s] = \E^{P_{\xi_s}}[\phi(\xi_{t-s}')].
\end{equation}
\end{theorem}
The set of times where the Markov property fails to hold at some point $x$
will be called the set of \emph{exceptional times} of $x$.

A very short outline of the proof (a complete version is given in~\cite{FlaRom08})
is the following:
\begin{itemize}
\item the sets $(\ms(x))_{x\in H}$ satisfy a \emph{set-valued} version of the Markov property,
\item given $\lambda>0$ and $f\in C_b(H)$, the set of maxima of the function
$$
P\mapsto\ris_{\lambda,f}(P):= \E^P\Bigl[\int_0^\infty\e^{-\lambda t}f(\xi_t)\,dt\Bigr]
$$
satisfies again the \emph{set-valued} version of the Markov property,
\item the proof is completed by iterating the above argument over a countable dense
set of $\lambda\in(0,\infty)$ and $f\in C_b(H)$.
\end{itemize}
In particular an arbitrary solution $P\in\ms(x)$ (for some $x\in H$) may not be an
element of a Markov solution (for a counterexample, see Proposition~\ref{p:peano_markov}).
Moreover, an arbitrary Markov solution may not be \emph{extremal},
that is, may not be obtained with the procedure outlined above (for a counterexample
see Proposition~\ref{p:peano_extremal}).

So far, the previous theorem ensures the existence of Markov solutions such that
the map $x\mapsto P_x$ is measurable (with respect to the appropriate $\sigma$-fields).
In order to improve the dependence with respect to the initial condition, the assumptions
on the covariance of the noise need to be strengthened. Assume that for some $\alpha_0>\frac16$
the operator $A^{\frac34+\alpha_0}\cov^{\frac12}$ is bounded invertible on H, with bounded inverse.
The additional regularity of noise trajectories allows to exploit the following fact: for
regular initial conditions there is a unique strong\footnote{The \emph{strong} has to be understood
in the PDE sense} solution up to a random time in which the derivatives of $u$ blow up.
The random time can be approximated monotonically by a sequence of stopping times
$$
\tau_x^\er = \inf\{t>0:\|u^\er_x(t)\|_{\W_{\alpha_0}}^2\geq R\}
$$
(see below for the definition of $\W_{\alpha_0}$), where $u^\er_x$ solves
\begin{equation}\label{e:nseR}
du_x^\er + \bigl(\nu Au_x^\er + \chi(\tfrac1R\|u_x^\er\|_{\W_{\alpha_0}}^2)B(u_x^\er,u_x^\er)\bigr)\,dt = \cov^\frac12\,dW,
\end{equation}
with initial condition $x$ and $\chi$ is a \emph{cut-off} function such that
$\chi(r)\equiv1$ for $r\leq1$ and $\chi(r)\equiv0$ for $r\geq2$.

Consider a Markov solution $(P_x)_{x\in H}$ and define for every $t\geq0$ the transition semigroup
$\semi_t:\bm(H)\to\bm(H)$ as
$$
\semi_t\varphi(x) = \E^{P_x}[\varphi(\xi_t)].
$$
\begin{theorem}[Theorem 5.11, \cite{FlaRom08}]\label{t:strongfeller}
Given a Markov solution $(P_x)_{x\in H}$, let $(\semi_t)_{t\geq0}$ be its transition semigroup.
Then for every $t>0$ and $\varphi\in\bm(H)$, $\semi_t\varphi\in C_b(\W_{\alpha_0})$.
\end{theorem}
The continuity in the above theorem is in the topology of $\W_{\alpha_0}=D(A^\theta)$,
where $\theta=\frac12(\alpha_0+1)$ if $\alpha_0<\frac12$ and
$\theta=\alpha_0+\frac14$ if $\alpha_0\geq\frac12$.
\begin{remark}[On regularity]\label{r:regularity}
Indeed, it is possible to improve the regularity result given above by relaxing the topology.
Under the same assumptions on the noise, one can replace in Theorem~\ref{t:strongfeller}
$\W_{\alpha_0}$ with $D(A^{1/4+\epsilon})$ (with arbitrary $\epsilon>0$), by using
parabolic regularisation. The main idea is that $D(A^{\frac14})$ is a critical space,
that is to say, it is the smallest space, in the hierarchy of hilbertian Sobolev spaces,
where it is possible to prove existence and uniqueness of local strong solutions for
the deterministic version of the problem. This extension is part of a work which
is currently in progress.
\end{remark}
\begin{remark}[On non-degeneracy]
The non-degeneracy assumption on the covariance $\cov$ can be slightly relaxed.
Assume for instance that $\cov$ is diagonal with respect to the Fourier basis
and assume that the range of $\cov$ spans all but a finite number of Fourier
modes. It is possible then to prove that any Markov solution is strong Feller
(this is a work in progress in collaboration with L. Xu). Unfortunately, at
least with the method presented here where the strong Feller property is
crucial, it does not seem to be possible to consider a noise highly
degenerate as, for example, in Hairer \& Mattingly~\cite{HaiMat06}.
\end{remark}
The above regularity result allows to analyse the long time behaviour
of any arbitrary Markov solution. The following theorem collects
Corollary 3.2 and Theorem 3.3 from~\cite{Rom08}.
\begin{theorem}\label{t:invariant}
Under the above assumptions on the covariance, every Markov solution
$(P_x)_{x\in H}$ to \eqref{e:nse} has a unique invariant measure $\mu$.
Moreover, there are $c>0$ nd $a>0$ such that
$$
\|\semi_t^*\delta_{x_0}-\mu\|_{\mathsf{TV}}\leq c\e^{-at},
$$
for every $x_0\in H$, where $\|\cdot\|_{\mathsf{TV}}$ is the
total variation norm.
\end{theorem}
It is worth remarking that the above theorem states that uniqueness of the
invariant measure holds among all measures which are invariant with respect
to the \emph{given} Markov solution. In different words, different
Markov solutions have different invariant measures. The following
result, which gathers Corollary $3.5$ and Theorem $4.6$ of~\cite{Rom08},
allows to compare invariant measures for different Markov solutions.
\begin{theorem}\label{t:equivalent}
All invariant measures associated to Markov solutions are mutually
equivalent. Moreover, if all Markov solutions share the same
invariant measure, then the martingale problem is well-posed.
\end{theorem}
\section{The martingale problem for the Navier-Stokes equations}\label{s:mp}
In this section we follow closely Appendix B of Cerrai~\cite{Cer01}.
Let $(P_x)_{x\in H}$ be a Markov solution and let $(\semi_t)_{t\geq0}$
be the associated transition semigroup. In the rest of the section
we will denote by $\W$ the space where the Markov solution is continuous,
without stating any dependence on the parameter $\alpha_0$.
In view of Remark~\ref{r:regularity}, there is no loss of generality
in doing so.
\begin{lemma}\label{l:continuous}
Let $\varphi\in\bm(H)$, then
\begin{itemize}
\item if $x\in\W$, then the map $t\mapsto\semi_t\varphi(x)$ is continuous for all $t\in[0,\infty)$,
\item if $x\in H$, then the map $t\mapsto\semi_t\varphi(x)$ is continuous for all $t\in(0,\infty)$.
\end{itemize}
\end{lemma}
\begin{proof}
If $x\in\W$, the statement follows from Lemma 6.6 of~\cite{FlaRom08}. If $x\in H$ and $t_0>0$,
choose $\delta>0$ such that $t_0-\delta>0$ and $t_0-\delta$ is not an \emph{exceptional time}
for $x$, then
$$
\semi_{t}\varphi(x)
 = \E^{P_x}[(\semi_{t-\delta}\varphi)(\xi_\delta)].
$$
Since by Lemma 3.7 of~\cite{Rom08} $P_x[\xi_\delta\in\W]=1$, by the first statement
of the lemma it follows that
$(\semi_{t-\delta}\varphi)(\xi_\delta)\to(\semi_{t_0-\delta}\varphi)(\xi_\delta)$
$P_x$-a.\ s.. The conclusion follows from Lebesgue theorem.
\end{proof}
Consider now $\lambda>0$ and define the operator $\ris_\lambda:\bm(H)\to\bm(H)$ as
$$
\ris_\lambda\varphi(x) = \int_0^\infty\e^{-\lambda t}\semi_t\varphi(x)\,dt.
$$
\begin{lemma}
For every $\lambda>0$, $\ris_\lambda$ is a bounded operator on $C_b(\W)$.
Moreover, the \emph{resolvent identity} holds. For every $\lambda_1$, $\lambda_2$,
$$
\ris_{\lambda_1} - \ris_{\lambda_2} = (\lambda_2 - \lambda_1)\ris_{\lambda_1}\ris_{\lambda_2}.
$$
\end{lemma}
\begin{proof}
Continuity of $\ris_\lambda\varphi$ follows from the strong Feller property and Lebesgue
theorem (since $\e^{-\lambda t}$ is integrable). Moreover, if $x\in\W$,
$$
|\ris_\lambda\varphi(x)|
\leq \int_0^\infty \e^{-\lambda t}|\semi_t\varphi(x)|\,dt
\leq \frac1\lambda\|\varphi\|_\infty.
$$
Next, we prove the resolvent identity. Fix $x\in\W$, then for a.\ e.\ $s$,
$\semi_{t+s}\varphi(x)=\semi_t\semi_s\varphi(x)$, and so
\begin{equation}\label{e:commute}
\begin{aligned}
\semi_t\ris_{\lambda_2}\varphi(x)
& = \E^{P_x}[\ris_{\lambda_2}\varphi(\xi_t)]
 = \E^{P_x}[\int_0^\infty\e^{-\lambda_2 s}\semi_s\varphi(\xi_t)\,ds]\\
& = \int_0^\infty\e^{-\lambda_2 s}\E^{P_x}[\semi_s\varphi(\xi_t)]\,ds
 = \int_0^\infty\e^{-\lambda_2 s}\semi_{t+s}\varphi(x)\,ds.
\end{aligned}
\end{equation}
In conclusion,
$$
\begin{aligned}
(\lambda_2-\lambda_1)\ris_{\lambda_1}\ris_{\lambda_2}\varphi(x)
& = (\lambda_2-\lambda_1)\int_0^\infty\e^{-\lambda_1 t}\int_0^\infty\e^{-\lambda_2 s}\semi_{t+s}\varphi(x)\,ds\,dt\\
& = (\lambda_2-\lambda_1)\int_0^\infty\e^{-\lambda_1 t}\int_t^\infty\e^{-\lambda_2 (r-t)}\semi_r\varphi(x)\,dr\,dt\\
& = \int_0^\infty\e^{-\lambda_2 r}(\e^{(\lambda_2-\lambda_1)r}-1)\semi_r\varphi(x)\,dr\\
& = \ris_{\lambda_1}\varphi(x) - \ris_{\lambda_2}\varphi(x),
\end{aligned}
$$
and the identity holds.
\end{proof}
We are finally able to prove existence of the generator.
\begin{theorem}\label{t:generator}
Given a Markov solution $(P_x)_{x\in H}$, there exists a unique closed
linear operator $\gen:D(\gen)\subset C_b(\W)\to C_b(\W)$ such that
for all $\lambda>0$ and $\varphi\in C_b(\W)$,
\begin{equation}\label{e:resolvent}
R_\lambda(\gen)\varphi(x) = \int_0^\infty\e^{-\lambda t}\semi_t\varphi(x)\,dt,
\end{equation}
where $R_\lambda(\gen)$ is the resolvent of $\gen$.
\end{theorem}
\begin{proof}
By the previous lemma, $(\ris_\lambda)_{\lambda>0}$ satisfies the resolvent
identity. Theorem $VIII.4.1$ of Yosida~\cite{Yos80}
ensures then that $(\ris_\lambda)_{\lambda>0}$ is the resolvent of a linear
operator $\gen$ if the kernel $N(\ris_\lambda)=\{0\}$. In such a case,
the domain $D(\gen)$ is equal to the range $R(\ris_\lambda)$, which is
independent of $\lambda$ by the resolvent identity.

We prove that $N(\ris_\lambda)=\{0\}$. Fix $\lambda_0>0$ and let $\varphi$
be such that $\ris_{\lambda_0}\varphi=0$. By the resolvent identity it
follows that $\ris_\lambda\varphi=0$ for all $\lambda>0$. By inverting
the Laplace transform, it follows that $\semi_t\varphi(x)=0$ for all
$x\in\W$ and almost every $t>0$ (hence all $t\geq0$ by Lemma~\ref{l:continuous}).
In particular, $\varphi=\semi_0\varphi=0$.
\end{proof}
\subsection{The martingale problem}

The computations of the previous section ensure that each Markov solution
has a generator. This allows to define the martingale problem associated
to this operator.
\begin{definition}[Martingale problem]
Let $\gen$ be the generator associated to some Markov solution and
provided by Theorem~\ref{t:generator} and let $x\in\W$. A probability
measure $P$ on $(\Omega,\B)$ is a solution to the martingale problem
associated to $\gen$ and starting at $x$ if
\begin{itemize}
\item $P[\xi_0=x]=1$,
\item for every $\varphi\in D(\gen)$, the process
      $$
      \mathcal{M}^\varphi_t = \varphi(\xi_t) - \int_0^t\gen\varphi(\xi_s)\,ds
      $$
      is a $P$-martingale with respect to the natural filtration $(\B_t)_{t\geq0}$.
\end{itemize}
\end{definition}
The aim of this section is to prove that each Markov solution is the
unique solution to the martingale problem associated to the corresponding
generator. With this aim in mind, we need the following lemma, which is
from Appendix~B of Cerrai~\cite{Cer01}. We give a short account of its
proof (which is essentially the same) because the assumptions under which
we work are slightly different.
\begin{lemma}\label{l:gen_is_sol}
For every $\varphi\in D(\gen)$ and $x\in\W$,
\begin{equation}\label{e:integform}
\semi_t\varphi(x) = \varphi(x) + \int_0^t\semi_s\gen\varphi(x)\,ds.
\end{equation}
In particular,
$$
\frac{d}{dt}\semi_t\varphi(x) = \semi_t\gen\varphi(x) = \gen\semi_t\varphi(x).
$$
\end{lemma}
\begin{proof}
By formula~\eqref{e:commute}, it follows that $\semi_t\ris_\lambda=\ris_\lambda\semi_t$.
Hence, $\semi_t(D(\gen))\subset D(\gen)$ since $D(\gen)=\ris_\lambda(C_b(\W))$ and so
$\gen\semi_t=\semi_t\gen$.

We prove~\eqref{e:integform}. Let $\varphi\in D(\gen)$, $x\in\W$ and $\lambda>0$, then
$$
\varphi(x)
= \ris_\lambda(\lambda I-\gen)\varphi(x)
=   \lambda\int_0^\infty\e^{-\lambda t}\semi_t\varphi(x)\,dt
  - \int_0^\infty\e^{-\lambda t}\semi_t\gen\varphi(x)\,dt
$$
and so by Fubini theorem,
$$
\begin{aligned}
\int_0^\infty\e^{-\lambda t}(\semi_t\varphi(x)-\varphi(x))\,dt
&= \int_0^\infty\frac1\lambda\e^{-\lambda t}\semi_t\gen\varphi(x)\,dt\\
&= \int_0^\infty\int_t^\infty\e^{-\lambda s}\semi_t\gen\varphi(x)\,ds\,dt\\
&= \int_0^\infty\e^{-\lambda s}\int_0^\infty\semi_t\gen\varphi(x)\,dt\,ds.
\end{aligned}
$$
By inverting the Laplace transform and using Lemma~\ref{l:continuous},
\eqref{e:integform} follows.
\end{proof}
\begin{theorem}
Let $(P_x)_{x\in\W}$ be a Markov solution and let $\gen$ be the associated
generator. Then the family $(P_x)_{x\in\W}$ is the unique solution to the
martingale problem associated to $\gen$.
\end{theorem}
\begin{proof}
Both proofs of existence and uniqueness are classical (see for example
Stroock \& Varadhan~\cite{StrVar79}), we give a proof for the interested
reader. First, we prove that $(P_x)_{x\in\W}$ is a solution to the martingale
problem. The Markov property~\eqref{e:markov} ensures that
$$
\E^{P_x}[\varphi(\xi_t)-\semi_{t-s}\varphi(\xi_s)|\B_s] = 0,
$$
while Lemma~\ref{l:gen_is_sol} implies that, $P_x$-a.\ s.,
$$
\varphi(\xi_s) = \semi_{t-s}\varphi(\xi_s) + \int_s^t\semi_{r-s}\gen\varphi(\xi_s)\,dr.
$$
Hence,
$$
\begin{aligned}
\E^{P_x}[\mathcal{M}_t^\varphi - \mathcal{M}_s^\varphi|\B_s]
&=\E^{P_x}\bigl[\varphi(\xi_t) - \varphi(\xi_s) - \int_s^t\gen\varphi(\xi_r)\,dr|\B_s\bigr]\\
&=\E^{P_x}\bigl[\varphi(\xi_t) - \semi_{t-s}\varphi(\xi_s) - \int_s^t\bigl(\gen\varphi(\xi_r) - \semi_{r-s}\gen\varphi(\xi_s)\bigr)|\B_s\bigr]\\
&=0.
\end{aligned}
$$
Next, we prove that $P_x$ is the unique solution. Let $P$ be a solution
to the martingale problem starting at $x$, let $\phi\in C_b(\W)$ and
set $\varphi = R_\lambda(\gen)\phi$. By definition of solution,
$$
\varphi(x) = \E^P[\varphi(\xi_t)-\int_0^t\gen\varphi(\xi_s)\,ds]
$$
and so by multiplying by $\lambda\e^{-\lambda t}$ and integrating by parts,
$$
\begin{aligned}
\varphi(x)
&=\E^P\bigl[\int_0^\infty\e^{-\lambda t}(\lambda I-\gen)\varphi(\xi_t)\,dt\bigr]\\
&=\E^P\bigl[\int_0^\infty\e^{-\lambda t}\phi(\xi_t)\,dt\bigr].
\end{aligned}
$$
By using~\eqref{e:resolvent} and inverting the Laplace transform, it follows
that $\E^{P}[\phi(\xi_t)]=\semi_t\phi(x)$. Since for a Markov process
uniqueness of one-dimensional distributions implies uniqueness of laws,
the theorem is proved.
\end{proof}
\begin{remark}
Da Prato and Debussche~\cite{DapDeb08} give a stricter definition
of solution to the martingale problem, due to the better knowledge
they have on their Markov solution, which is obtained via Galerkin
approximations (see~\cite{DapDeb03}).
\end{remark}
\subsection{What can we say of the generator}\label{ss:generator}

So far, we have proved that any Markov solution $(P_x)_{x\in H}$ is the
unique solution to the martingale problem associated to the generator
of the transition semigroup. On the other hand, the formal expression
of the generator associated to~\eqref{e:abstractnse} is
$$
\gen^\star\varphi(x)
= \frac12\Tr[\cov D^2\varphi](x) - \scal{\nu Ax + B(x,x), D\varphi(x)}.
$$
In this section we shall try to understand (although without success)
if there is any relation between $\gen^\star$ and the generator $\gen$
of an arbitrary Markov solution.

To this aim, fix a Markov solution $(P_x)_{x\in H}$ and let $\gen$
be the associated generator. It is useful to notice that the generator
$\gen$ can be characterised (see Da Prato \& Debussche~\cite{DapDeb08})
in the following way,
$$
D(\gen) = \bigl\{\varphi\in C_b(\W):
	\lim_{\eps\to0}\frac{\semi_\eps\varphi(x)-\varphi(x)}{\eps}\text{ exists }\forall x\in\W\text{ and is in }C_b(\W)\bigr\}
$$
and
$$
\gen\varphi(x) = \lim_{\eps\to0}\frac{\semi_\eps\varphi(x)-\varphi(x)}{\eps}.
$$
Let $(\semi^\er_t)_{t\geq0}$ the Markov semigroup associated to the
cut-off problem~\eqref{e:nseR} and let $\gen^\er$ be the corresponding generator.
\begin{lemma}
Given $R\geq1$,
\begin{itemize}
\item if $\phi\in D(\gen^\er)$, then for every $|x|_\W^2< R$, $\gen^\er\phi(x)=\gen^\star\phi(x)$,
\item if $\phi\in D(\gen^\er)$, then for every $|x|_\W^2< R$,
      $$
      \lim_{t\to0}\tfrac1{t}(\semi_t\phi(x)-\phi(x))=\gen^\er\phi(x),
      $$
\item if $\phi\in D(\gen)$, then for every $|x|_\W^2< R$,
      $$
      \lim_{t\to0}\tfrac1{t}(\semi_t\phi(x)-\phi(x))=\gen\phi(x),
      $$
\end{itemize}
\end{lemma}
\begin{proof}
The first property is easy. The second and third property follow from
\begin{align}\label{e:equal_gen}
\Bigl|\frac{\semi_t^\er\phi(x)-\phi(x)}{t} - \frac{\semi_t\phi(x)-\phi(x)}{t}\Bigr|
& =   \Bigl|\frac{\semi_t^\er\phi(x) - \semi_t\phi(x)}{t}\Bigr|\notag\\
&\leq \tfrac{2}{t}\|\phi\|_\infty P_x^\er[\tau_x^\er<t]\\
&\leq \tfrac{c_1}{t}\|\phi\|_\infty \e^{-c_2\frac{R^2}{t}},\notag
\end{align}
for $t\leq cR^{-\gamma}$ (for some $c>0$ and $\gamma>0$), where the first inequality
follows from Lemma 5.9 of~\cite{FlaRom08} (see also part $1$ of Theorem~A.1 in~\cite{Rom08})
and the second inequality follows from Proposition 11 of~\cite{FlaRom07} (see also part
$2$ of Theorem~A.1 in~\cite{Rom08}).
\end{proof}
Based on this lemma, the following proposition gives a (almost elementary)
condition for the generator $\gen$ to be equal to the formal expression $\gen^\star$.
\begin{proposition}
Let
$$
\mathcal{E}
=\{\phi\in D(\gen): \text{there is }R_n\uparrow\infty\text{ s.t. }\phi\in\bigcap_{n\in\N}D(\gen^\ern)\text{ and }\sup_{n\in\N}\|\gen^\ern\phi\|_\infty\}.
$$
Then $\gen\phi = \gen^\star\phi$ for every $\phi\in\mathcal{E}$.
\end{proposition}
\begin{proof}
The property follows from inequality~\eqref{e:equal_gen}, since $\tfrac{1}{t}(\semi_t^\ern\phi(x)-\phi(x))$ 
is bounded because
$$
\Bigl|\frac{\semi_t^\ern\phi(x)-\phi(x)}{t}\Bigr|
=    \Bigl|\int_0^t\semi_s^\ern\gen^\ern\phi(x)\,ds\Bigr|
\leq \sup_{n\in\N}\|\gen^\ern\phi\|_\infty
$$
and $\tfrac{1}{t}(\semi_t\phi(x)-\phi(x))$ is bounded by the alternative
characterisation of $\gen$ given above.
\end{proof}
\begin{remark}\label{r:bad_tau}
No better conclusion can be drawn with such generality (see Da Prato
\& Debussche~\cite{DapDeb08} for some related results). The argument
missing in this analysis is, essentially, a better estimate of tails
of the stopping time $\tau_x^\er$, which is used in formula~\eqref{e:equal_gen}
to estimate the distance from the generator to the cut-off problem.
\end{remark}
\section{Some elementary examples of Markov solutions}\label{s:examples}

In this last section we present some elementary examples from the theory
of (deterministic and stochastic) differential equations. We wish to
compare these with all results on the stochastic Navier-Stokes equations
given in the previous sections.

The first example is a revisitation of a classical example of non-uniqueness
in ordinary differential equations, where it is easy to characterise
all Markov solutions (compare Proposition~\ref{p:peano_extremal} with
Theorem~\ref{t:equivalent}).

The second example is taken from a paper by Girsanov~\cite{Gir62}, wihere
all Markov solutions are Feller and it is possible to list the
generators of all such solutions (compare with Section~\ref{ss:generator}).

The last example has been presented by Stroock \& Yor~\cite{StrYor80}
and its main interest is that there are two (strong) Markov solutions
which are both strong Feller (compare with Theorem~\ref{t:strongfeller}).
\subsection{An example from elementary calculus}\label{ss:peano}
Consider the following differential equation
\begin{equation}\label{e:peano}
\dot X = -X + \sqrt{X},
\end{equation}
with initial condition $X(0)=x\in[0,1]$. The problem has
a unique solution $X_x(\cdot)$ for $x\neq0$ and the family
of solutions
$$
\{X_a^\star=X^\star((t-a)\vee0):a\geq0\}
$$
for $x=0$, where $X^\star$ is the unique solution starting
at $0$ such that $X^\star(t)>0$ for all $t>0$.

If $\msP(x)$ denotes the set of all solutions to~\eqref{e:peano} starting at $x$,
then $\msP(x)=\{\delta_{X_x}\}$ for $x\in(0,1]$, where $\delta_{X_x}$
is the Dirac measure on $C([0,\infty);\Erre)$ concentrated on $X_x$.

\begin{figure}
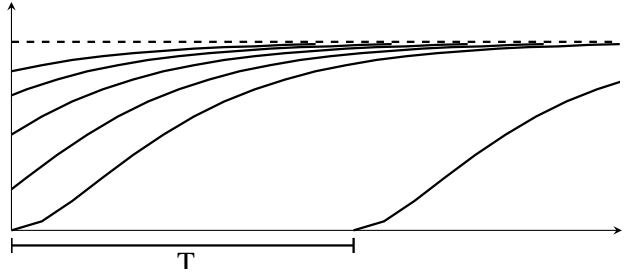
\label{f:peano}
\begin{pgfpicture}{0mm}{0mm}{80mm}{35mm}
  \pgfsetxvec{\pgfpoint{1mm}{0mm}}
  \pgfsetyvec{\pgfpoint{0mm}{1mm}}
  \pgfsetendarrow{\pgfarrowsingle}
  \pgfline{\pgfxy(0,5)}{\pgfxy(80,5)}
  \pgfline{\pgfxy(0,5)}{\pgfxy(0,35)}
  \pgfsetendarrow{}
  \pgfpathmoveto{\pgfxy(0,0)}
  \pgfpathlineto{\pgfxy(0,35)}
  \pgfpathlineto{\pgfxy(80,35)}
  \pgfpathlineto{\pgfxy(80,0)}
  \pgfclip
  \pgfsetlinewidth{0.3mm}
  \curva{0,5}
  \curva{-10,5}
  \curva{-20,5}
  \curva{-30,5}
  \curva{-40,5}
  \curva{45,5}
  \pgfline{\pgfxy(0,2)}{\pgfxy(0,4)}
  \pgfline{\pgfxy(45,2)}{\pgfxy(45,4)}
  \pgfline{\pgfxy(0,3)}{\pgfxy(45,3)}
  \pgfputat{\pgfxy(23,2)}{\pgfbox[center,top]{\small T}}
  \pgfsetdash{{1mm}}{0pt}
  \pgfline{\pgfxy(0,30)}{\pgfxy(80,30)}
\end{pgfpicture}
\caption{Solutions of the Peano example}
\end{figure}

If $x=0$, a solution starts at $0$ and stays for an arbitrary time,
then follows the solution $X^\star$ (suitably translated). So the
\emph{departing} time from $0$ can be interpreted as a random
time whose law can be arbitrary (see fig.~\ref{f:peano}). 
\begin{lemma}
The set of solutions starting at $x=0$ is given by
$$
\msP(0) = \Bigl\{\int\delta_{X^\star_a}\,\mu(da):\mu\text{ is a probability measure on }[0,\infty]\Bigr\}
$$
\end{lemma}
In conclusion, any selection is completely described by a single random
variable on $[0,\infty)$ (or, equivalently, by a single measure on $[0,\infty)$).
Given a probability measure $\nu$ on $[0,\infty)$, define
$$
P^\nu_x=
\begin{cases}
\delta_{X_x}&\qquad x\in(0,1],\\
\int\delta_{X_a^\star}\nu(da)&\qquad x=0,
\end{cases}
$$
then $(P^\nu_x)_{x\in[0,1]}$ is a measurable selection, and any selection
corresponds to one of them for some $\nu$.
\begin{proposition}\label{p:peano_markov}
A measurable selection $(P^\nu_x)_{x\in[0,1]}$ is Markov if and only if $\nu$
is the distribution of an exponential random variable (including the degenerate
cases of infinite or zero rate, where $\nu=\delta_0$ or $\nu=\delta_\infty$).
\end{proposition}
\begin{proof}
It is easy to see that the Markov property holds if $x\in(0,1]$ whatever is $\nu$.
Indeed, $P_x^\nu$-a.\ s.,
$$
\E^{P^\nu_x}[f(\xi_{t+s})]
=f(X_x(t+s))
=f(X_{X_x(t)}(s))
=\E^{P^\nu_x}[f(X_{\xi_t}(s))]
=\E^{P^\nu_x}[\E^{P^\nu_{\xi_t}}[f(\xi_s')]].
$$
We next see which condition we get if $x=0$. On one side,
\begin{align*}
\E^{P^\nu_0}[f(\xi_{t+s})]
&=\int f(X_a^\star(t+s))\nu(da)\\
&=\int_{[0,t)}f(X^\star(t+s-a))\nu(da) + \int_{[t,+\infty]}f(X^\star((t+s-a)\vee0))\nu(da)\\
&=\framebox{1} + \int f(X_b^\star(s))(\theta_t\nu)(db),
\end{align*}
where $\theta_t:[t,\infty]\to[0,\infty]$ is defined as $\theta_t(s)=s-t$.
On the other side,
\begin{align*}
\E^{P^\nu_0}[\E^{P^\nu_{\xi_t}}[f(\xi_s')]]
&=\int_{[0,t)}\E^{P^\nu_{X^\star(t-a)}}[f(\xi_s)]\nu(da) + \int_{[t,+\infty]}\E^{P^\nu_0}[f(\xi_s)]\nu(da)\\
&=\int_{[0,t)}f(X^\star(t+s-a))\nu(da) + \nu([t,\infty])\int f(X_b^\star(s))\nu(db)\\
&=\framebox{1} + \nu([t,\infty])\int f(X_b^\star(s))\nu(db).
\end{align*}
In conclusion,
$$
\int f(X_b^\star(s))(\theta_t\nu)(db)
=\nu([t,\infty])\int f(X_b^\star(s))\nu(db).
$$
Moreover, by splitting the integrals in the formula above on $[0,s)$ and $[s,\infty]$,
\begin{multline*}
f(0)\phi(s+t)+\int_{[0,s)}f(X^\star(s-b))(\theta_t\nu)(db)=\\
=f(0)\phi(s)\phi(t) + \phi(t)\int_{[0,s)}f(X^\star(s-b))\nu(db),
\end{multline*}
where $\phi(r)=\nu([r,\infty])$.

A further simplification can be achieved since $b\in(0,s]\to X^\star(s-b)\in(0,X^\star(s)]$
is invertible (with inverse $g$ and $g(0)=s$), so that if $f=F\circ g$, then we finally get
$$
F(s)\phi(s+t)+\int_{[0,s)}F(b)(\theta_t\nu)(db)
=F(s)\phi(s)\phi(t) + \phi(t)\int_{[0,s)}F(b)\nu(db).
$$
This implies that $\phi(t)\nu=\theta_t\nu$ and
$$
\phi(s+t)=\phi(s)\phi(t)
$$
and $\nu$ is the law of an exponential random variable.
\end{proof}
For every $a\in[0,\infty]$ we denote by $(P^a_x)_{x\in[0,1]}$ the Markov selection
of rate $a$. We shall call \emph{extremal} all those Markov solutions that can be
obtained by the selection procedure outlined in Section~\ref{s:markov}.
\begin{proposition}\label{p:peano_extremal}
The extremal selections are those corresponding to $a=0$ and $a=\infty$.
\end{proposition}
\begin{proof}
Given $\lambda>0$ and a function $f$, a straightforward computation gives
\begin{align*}
J_{\lambda,f}(P_0^a)-J_{\lambda,f}(P_0^b)
&=\frac{b-a}{(\lambda+a)(\lambda+b)}[\lambda J_{\lambda,f}(P_0^0)-f(0)]\\
&=\frac{\lambda(b-a)}{(\lambda+a)(\lambda+b)}[J_{\lambda,f}(P_0^0)-J_{\lambda,f}(P_0^\infty)],
\end{align*}
and with this formula the conclusion is obvious.
\end{proof}
As it regards invariant measures, we notice that $(P_x^a)_{x\in [0,1]}$
has a unique invariant measure (which is $\delta_1$) if and only if
$a<\infty$. Notice that all selections having $\delta_1$ as their
unique invariant measure coincide $\delta_1$--almost surely.

If $a=\infty$, there are infinitely many invariant measures
(the convex hull of $\delta_0$ and $\delta_1$). As there is
no noise in this example, in general we cannot expect the
invariant measures to be equivalent (compare with
Theorem~\ref{t:equivalent}).
\subsection{An example of non-uniqueness from Girsanov}\label{ss:girsanov}
In his paper~\cite{Gir62}, Girsanov is able to classify
the generators of all diffusions which solve the following
stochastic differential equation,
\begin{equation}\label{e:girsanov}
dX_t = \sigma_\alpha(X_t)\,dW_t,
\end{equation}
where, for any $\alpha\in(0,\tfrac12)$, $\sigma_\alpha$ is the function
$$
\sigma_\alpha(x)=\frac{|x|^\alpha}{1+|x|^\alpha}.
$$
Engelbert \& Schmidt~\cite{EngSch85} give a characterisation
for existence and uniqueness of one dimensional \emph{SDE}s
as the one under examination. Their Theorem $2.2$ ensures that
there is at least one solution for each initial condition, while
their Theorem $3.2$ implies that the problem has no unique
solution.
\begin{remark}
The same conclusions hold for a generic function $\sigma$
such that $\sigma^{-2}$ is locally integrable and the
set of zeroes $\{x: \sigma^2(x)=0\}$ is not empty.
See also Example $4.1$ of Stroock \& Yor~\cite{StrYor80}.
\end{remark}
In the rest of this section we give a twofold description
of Markov solutions to problem~\eqref{e:girsanov}, in
terms of the generator and in terms of the process.
\subsubsection{The generators}

Girsanov~\cite{Gir62} shows that each of the Markov solutions
has its own generator $\genG$ with domain $D(\genG)$. All functions
$C_b^2(\Erre)$ are in $D(\genG)$ and for $x\neq0$,
$$
\genG u(x) = \sigma_\alpha(x)^2 u''(x).
$$
If the solution corresponds to the point $0$ to be \emph{absorbing}
(i.\ e., the solutions stays in $0$ once it hits it),
then
$$
\genG_\infty u(0) = 0.
$$
In the \emph{non-absorbing} case, the generators can be parametrised
by $c\geq0$. If $c>0$, the domain $D(\genG_c)$ contains all $C_b^2(\Erre\backslash\{0\})$
such that the left and right derivatives exist in $0$,
$$
\genG_c u(0) = \frac{1}{c}\bigl(u'(0+)-u'(0-)\bigr),
$$
and $\genG_c u$ is continuous on $\Erre$. In the case $c=0$,
which corresponds to a process which spends no time in $0$,
$$
\genG_0 u(0) = \lim_{x\to0} \sigma_\alpha(x)^2 u''(x).
$$
The meaning of the parameter $c$ will be clarified in the next
section, where we shall give an explicit construction of Markov
processes solving the problem (see also McKean~\cite[Section 3.10b]{McK69}).
\subsubsection{Description of solutions via time-changes and delays}

Following Theorem 2.2 of Engelbert \& Schmidt~\cite{EngSch85},
we start by the construction of a process corresponding to
$c=0$. Define the (strictly increasing) process
$$
S_t^x = \int_0^t\frac1{\sigma_\alpha(x+W_s)^2}\,ds,
$$
(the integrability of $\sigma_\alpha^{-2}$ ensures that $S_t^x<\infty$
for all $t\geq0$, $\Pb$-a.s.), and denote by $T_t^x$ the inverse of $S^x$.
The process $(T_t^x)_{t\geq0}$ is again strictly increasng and $T_\infty^x=\infty$.
By Proposition $5.1$ of Stroock \& Yor~\cite{StrYor80}, the process $\Xz_t^x = x + W_{T_t^x}$
is a solution to~\eqref{e:girsanov}. Moreover, by Theorem $5.4$ of \cite{EngSch85},
it is the only solution such that
\begin{equation}\label{e:notimeinzero}
\E\Bigl[\int_0^{+\infty}\uno_{\{0\}}(X_t)\,dt\Bigr] = 0.
\end{equation}
The above condition \eqref{e:notimeinzero} means that $\Xz^x$ spends no
time in $0$ and ensures that $(\Xz_t^x)_{x\in\Erre}$ is a Markov process.

Any other solution can be obtained by \emph{delaying} $\Xz$ (Theorem $5.5$
of~\cite{EngSch85}). Indeed, a \emph{time-delay} for $\Xz$ is any adapted
increasing right-continuous process $(D_t)_{t\geq0}$ such that
$$
D_t = \int_0^t\uno_{\{0\}}(\Xz_s)\,dD_s,
\qquad t\geq0,\quad\Pb-\text{q.c.}.
$$
If $E_t$ is the inverse of $t\mapsto t+D_t$, then the process $Y_t=\Xz_{E_t}$,
adapted to $\F_t=\B_{E_t}$ is again a solution (Theorem $4.3$ of~\cite{EngSch85}).
\begin{figure}[ht]
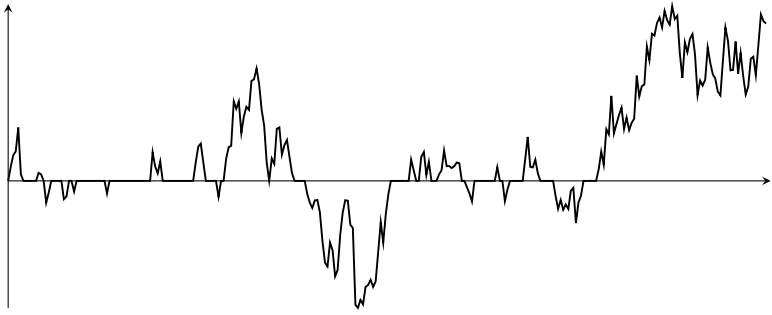

\begin{pgfpicture}{0mm}{0mm}{100mm}{40mm}
  \pgfsetxvec{\pgfpoint{1mm}{0mm}}
  \pgfsetyvec{\pgfpoint{0mm}{1mm}}
  \pgfsetendarrow{\pgfarrowsingle}
  \pgfline{\pgfxy(0,0)}{\pgfxy(0,40)}
  \pgfline{\pgfxy(0,16.84)}{\pgfxy(100,16.84)}
  \pgfsetendarrow{}
  \pgfsetlinewidth{0.25mm}
  \pgfpathmoveto{\pgfxy( 0.00, 16.84)}\pgfpathlineto{\pgfxy( 0.33, 18.67)}\pgfpathlineto{\pgfxy( 0.67, 20.24)}\pgfpathlineto{\pgfxy( 1.00, 20.76)}
  \pgfpathlineto{\pgfxy( 1.33, 23.91)}\pgfpathlineto{\pgfxy( 1.67, 17.72)}\pgfpathlineto{\pgfxy( 2.00, 16.84)}\pgfpathlineto{\pgfxy( 3.67, 16.84)}
  \pgfpathlineto{\pgfxy( 4.00, 17.90)}\pgfpathlineto{\pgfxy( 4.33, 17.71)}\pgfpathlineto{\pgfxy( 4.67, 16.84)}\pgfpathlineto{\pgfxy( 5.00, 13.94)}
  \pgfpathlineto{\pgfxy( 5.33, 15.26)}\pgfpathlineto{\pgfxy( 5.67, 16.84)}\pgfpathlineto{\pgfxy( 7.00, 16.84)}\pgfpathlineto{\pgfxy( 7.33, 14.41)}
  \pgfpathlineto{\pgfxy( 7.67, 14.77)}\pgfpathlineto{\pgfxy( 8.00, 16.84)}\pgfpathlineto{\pgfxy( 8.33, 16.84)}\pgfpathlineto{\pgfxy( 8.67, 15.42)}
  \pgfpathlineto{\pgfxy( 9.00, 16.84)}\pgfpathlineto{\pgfxy(12.67, 16.84)}\pgfpathlineto{\pgfxy(13.00, 15.16)}\pgfpathlineto{\pgfxy(13.33, 16.84)}
  \pgfpathlineto{\pgfxy(18.67, 16.84)}\pgfpathlineto{\pgfxy(19.00, 20.63)}\pgfpathlineto{\pgfxy(19.33, 18.78)}\pgfpathlineto{\pgfxy(19.67, 17.84)}
  \pgfpathlineto{\pgfxy(20.00, 19.42)}\pgfpathlineto{\pgfxy(20.33, 16.84)}\pgfpathlineto{\pgfxy(24.33, 16.84)}\pgfpathlineto{\pgfxy(24.67, 19.30)}
  \pgfpathlineto{\pgfxy(25.00, 21.37)}\pgfpathlineto{\pgfxy(25.33, 21.79)}\pgfpathlineto{\pgfxy(25.67, 19.26)}\pgfpathlineto{\pgfxy(26.00, 16.84)}
  \pgfpathlineto{\pgfxy(27.33, 16.84)}\pgfpathlineto{\pgfxy(27.67, 14.73)}\pgfpathlineto{\pgfxy(28.00, 16.84)}\pgfpathlineto{\pgfxy(28.33, 16.84)}
  \pgfpathlineto{\pgfxy(28.67, 19.80)}\pgfpathlineto{\pgfxy(29.00, 21.30)}\pgfpathlineto{\pgfxy(29.33, 21.49)}\pgfpathlineto{\pgfxy(29.67, 27.38)}
  \pgfpathlineto{\pgfxy(30.00, 26.40)}\pgfpathlineto{\pgfxy(30.33, 27.37)}\pgfpathlineto{\pgfxy(30.67, 23.19)}\pgfpathlineto{\pgfxy(31.00, 25.39)}
  \pgfpathlineto{\pgfxy(31.33, 26.66)}\pgfpathlineto{\pgfxy(31.67, 26.25)}\pgfpathlineto{\pgfxy(32.00, 30.11)}\pgfpathlineto{\pgfxy(32.33, 30.30)}
  \pgfpathlineto{\pgfxy(32.67, 31.75)}\pgfpathlineto{\pgfxy(33.00, 29.73)}\pgfpathlineto{\pgfxy(33.33, 26.33)}\pgfpathlineto{\pgfxy(33.67, 24.13)}
  \pgfpathlineto{\pgfxy(34.00, 19.28)}\pgfpathlineto{\pgfxy(34.33, 16.84)}\pgfpathlineto{\pgfxy(34.67, 19.81)}\pgfpathlineto{\pgfxy(35.00, 19.10)}
  \pgfpathlineto{\pgfxy(35.33, 23.74)}\pgfpathlineto{\pgfxy(35.67, 23.93)}\pgfpathlineto{\pgfxy(36.00, 20.32)}\pgfpathlineto{\pgfxy(36.33, 21.57)}
  \pgfpathlineto{\pgfxy(36.67, 22.27)}\pgfpathlineto{\pgfxy(37.00, 20.00)}\pgfpathlineto{\pgfxy(37.33, 17.87)}\pgfpathlineto{\pgfxy(37.67, 16.84)}
  \pgfpathlineto{\pgfxy(39.00, 16.84)}\pgfpathlineto{\pgfxy(39.33, 15.13)}\pgfpathlineto{\pgfxy(39.67, 13.95)}\pgfpathlineto{\pgfxy(40.00, 13.24)}
  \pgfpathlineto{\pgfxy(40.33, 14.25)}\pgfpathlineto{\pgfxy(40.67, 14.31)}\pgfpathlineto{\pgfxy(41.00, 12.64)}\pgfpathlineto{\pgfxy(41.33,  8.89)}
  \pgfpathlineto{\pgfxy(41.67,  6.00)}\pgfpathlineto{\pgfxy(42.00,  5.49)}\pgfpathlineto{\pgfxy(42.33,  8.67)}\pgfpathlineto{\pgfxy(42.67,  7.62)}
  \pgfpathlineto{\pgfxy(43.00,  4.18)}\pgfpathlineto{\pgfxy(43.33,  5.03)}\pgfpathlineto{\pgfxy(43.67,  9.67)}\pgfpathlineto{\pgfxy(44.00, 12.64)}
  \pgfpathlineto{\pgfxy(44.33, 14.28)}\pgfpathlineto{\pgfxy(44.67, 14.21)}\pgfpathlineto{\pgfxy(45.00, 11.02)}\pgfpathlineto{\pgfxy(45.33, 10.61)}
  \pgfpathlineto{\pgfxy(45.67,  0.41)}\pgfpathlineto{\pgfxy(46.00,  0.00)}\pgfpathlineto{\pgfxy(46.33,  1.05)}\pgfpathlineto{\pgfxy(46.67,  0.46)}
  \pgfpathlineto{\pgfxy(47.00,  2.77)}\pgfpathlineto{\pgfxy(47.33,  3.01)}\pgfpathlineto{\pgfxy(47.67,  3.72)}\pgfpathlineto{\pgfxy(48.00,  2.77)}
  \pgfpathlineto{\pgfxy(48.33,  3.50)}\pgfpathlineto{\pgfxy(48.67,  7.49)}\pgfpathlineto{\pgfxy(49.00, 11.36)}\pgfpathlineto{\pgfxy(49.33,  8.66)}
  \pgfpathlineto{\pgfxy(49.67, 12.49)}\pgfpathlineto{\pgfxy(50.00, 15.06)}\pgfpathlineto{\pgfxy(50.33, 16.84)}\pgfpathlineto{\pgfxy(52.67, 16.84)}
  \pgfpathlineto{\pgfxy(53.00, 19.66)}\pgfpathlineto{\pgfxy(53.33, 18.26)}\pgfpathlineto{\pgfxy(53.67, 16.84)}\pgfpathlineto{\pgfxy(54.00, 16.84)}
  \pgfpathlineto{\pgfxy(54.33, 20.02)}\pgfpathlineto{\pgfxy(54.67, 20.67)}\pgfpathlineto{\pgfxy(55.00, 17.63)}\pgfpathlineto{\pgfxy(55.33, 19.36)}
  \pgfpathlineto{\pgfxy(55.67, 16.84)}\pgfpathlineto{\pgfxy(56.33, 16.84)}\pgfpathlineto{\pgfxy(56.67, 17.70)}\pgfpathlineto{\pgfxy(57.00, 18.26)}
  \pgfpathlineto{\pgfxy(57.33, 20.87)}\pgfpathlineto{\pgfxy(57.67, 18.79)}\pgfpathlineto{\pgfxy(58.00, 18.81)}\pgfpathlineto{\pgfxy(58.33, 18.54)}
  \pgfpathlineto{\pgfxy(58.67, 18.78)}\pgfpathlineto{\pgfxy(59.00, 19.28)}\pgfpathlineto{\pgfxy(59.33, 19.19)}\pgfpathlineto{\pgfxy(59.67, 16.84)}
  \pgfpathlineto{\pgfxy(60.00, 16.84)}\pgfpathlineto{\pgfxy(60.33, 16.02)}\pgfpathlineto{\pgfxy(60.67, 15.17)}\pgfpathlineto{\pgfxy(61.00, 14.11)}
  \pgfpathlineto{\pgfxy(61.33, 16.84)}\pgfpathlineto{\pgfxy(64.00, 16.84)}\pgfpathlineto{\pgfxy(64.33, 18.66)}\pgfpathlineto{\pgfxy(64.67, 16.84)}
  \pgfpathlineto{\pgfxy(65.00, 16.84)}\pgfpathlineto{\pgfxy(65.33, 14.15)}\pgfpathlineto{\pgfxy(65.67, 15.72)}\pgfpathlineto{\pgfxy(66.00, 16.84)}
  \pgfpathlineto{\pgfxy(67.67, 16.84)}\pgfpathlineto{\pgfxy(68.00, 19.73)}\pgfpathlineto{\pgfxy(68.33, 22.67)}\pgfpathlineto{\pgfxy(68.67, 18.69)}
  \pgfpathlineto{\pgfxy(69.00, 18.62)}\pgfpathlineto{\pgfxy(69.33, 19.65)}\pgfpathlineto{\pgfxy(69.67, 17.80)}\pgfpathlineto{\pgfxy(70.00, 16.84)}
  \pgfpathlineto{\pgfxy(71.67, 16.84)}\pgfpathlineto{\pgfxy(72.00, 14.81)}\pgfpathlineto{\pgfxy(72.33, 13.13)}\pgfpathlineto{\pgfxy(72.67, 14.31)}
  \pgfpathlineto{\pgfxy(73.00, 13.01)}\pgfpathlineto{\pgfxy(73.33, 13.71)}\pgfpathlineto{\pgfxy(73.67, 13.12)}\pgfpathlineto{\pgfxy(74.00, 15.52)}
  \pgfpathlineto{\pgfxy(74.33, 15.92)}\pgfpathlineto{\pgfxy(74.67, 11.24)}\pgfpathlineto{\pgfxy(75.00, 13.95)}\pgfpathlineto{\pgfxy(75.33, 14.85)}
  \pgfpathlineto{\pgfxy(75.67, 16.84)}\pgfpathlineto{\pgfxy(77.33, 16.84)}\pgfpathlineto{\pgfxy(77.67, 18.42)}\pgfpathlineto{\pgfxy(78.00, 20.71)}
  \pgfpathlineto{\pgfxy(78.33, 18.99)}\pgfpathlineto{\pgfxy(78.67, 23.69)}\pgfpathlineto{\pgfxy(79.00, 23.01)}\pgfpathlineto{\pgfxy(79.33, 28.11)}
  \pgfpathlineto{\pgfxy(79.67, 23.15)}\pgfpathlineto{\pgfxy(80.00, 24.38)}\pgfpathlineto{\pgfxy(80.33, 25.55)}\pgfpathlineto{\pgfxy(80.67, 26.55)}
  \pgfpathlineto{\pgfxy(81.00, 23.70)}\pgfpathlineto{\pgfxy(81.33, 25.23)}\pgfpathlineto{\pgfxy(81.67, 23.56)}\pgfpathlineto{\pgfxy(82.00, 24.54)}
  \pgfpathlineto{\pgfxy(82.33, 25.08)}\pgfpathlineto{\pgfxy(82.67, 30.79)}\pgfpathlineto{\pgfxy(83.00, 28.04)}\pgfpathlineto{\pgfxy(83.33, 29.39)}
  \pgfpathlineto{\pgfxy(83.67, 29.62)}\pgfpathlineto{\pgfxy(84.00, 34.66)}\pgfpathlineto{\pgfxy(84.33, 32.85)}\pgfpathlineto{\pgfxy(84.67, 36.35)}
  \pgfpathlineto{\pgfxy(85.00, 36.11)}\pgfpathlineto{\pgfxy(85.33, 37.77)}\pgfpathlineto{\pgfxy(85.67, 38.50)}\pgfpathlineto{\pgfxy(86.00, 37.19)}
  \pgfpathlineto{\pgfxy(86.33, 39.36)}\pgfpathlineto{\pgfxy(86.67, 38.15)}\pgfpathlineto{\pgfxy(87.00, 37.55)}\pgfpathlineto{\pgfxy(87.33, 40.00)}
  \pgfpathlineto{\pgfxy(87.67, 38.30)}\pgfpathlineto{\pgfxy(88.00, 38.73)}\pgfpathlineto{\pgfxy(88.33, 33.80)}\pgfpathlineto{\pgfxy(88.67, 30.49)}
  \pgfpathlineto{\pgfxy(89.00, 35.20)}\pgfpathlineto{\pgfxy(89.33, 33.90)}\pgfpathlineto{\pgfxy(89.67, 35.68)}\pgfpathlineto{\pgfxy(90.00, 36.30)}
  \pgfpathlineto{\pgfxy(90.33, 33.62)}\pgfpathlineto{\pgfxy(90.67, 28.22)}\pgfpathlineto{\pgfxy(91.00, 30.14)}\pgfpathlineto{\pgfxy(91.33, 29.47)}
  \pgfpathlineto{\pgfxy(91.67, 30.24)}\pgfpathlineto{\pgfxy(92.00, 34.53)}\pgfpathlineto{\pgfxy(92.33, 32.55)}\pgfpathlineto{\pgfxy(92.67, 30.98)}
  \pgfpathlineto{\pgfxy(93.00, 30.39)}\pgfpathlineto{\pgfxy(93.33, 28.68)}\pgfpathlineto{\pgfxy(93.67, 28.17)}\pgfpathlineto{\pgfxy(94.00, 32.75)}
  \pgfpathlineto{\pgfxy(94.33, 37.18)}\pgfpathlineto{\pgfxy(94.67, 35.40)}\pgfpathlineto{\pgfxy(95.00, 31.50)}\pgfpathlineto{\pgfxy(95.33, 31.55)}
  \pgfpathlineto{\pgfxy(95.67, 35.32)}\pgfpathlineto{\pgfxy(96.00, 31.06)}\pgfpathlineto{\pgfxy(96.33, 33.80)}\pgfpathlineto{\pgfxy(96.67, 30.73)}
  \pgfpathlineto{\pgfxy(97.00, 28.29)}\pgfpathlineto{\pgfxy(97.33, 29.32)}\pgfpathlineto{\pgfxy(97.67, 33.04)}\pgfpathlineto{\pgfxy(98.00, 33.31)}
  \pgfpathlineto{\pgfxy(98.33, 30.85)}\pgfpathlineto{\pgfxy(98.67, 34.88)}\pgfpathlineto{\pgfxy(99.00, 38.91)}\pgfpathlineto{\pgfxy(99.33, 38.04)}
  \pgfpathlineto{\pgfxy(99.67, 37.69)}
  \pgfusepath{stroke}
\end{pgfpicture}
\caption{The solution is delayed whenever it hits $0$}
\end{figure}

In particular, if $\tau_0=\inf\{t\geq0:\Xz_t=0\}$ and $D_t=0$ for $t<\tau_0$,
and $+\infty$ otherwise, then the delayed process is the process stopped at $0$
(which corresponds to the generator $\genG_\infty$).

\begin{figure}[ht]
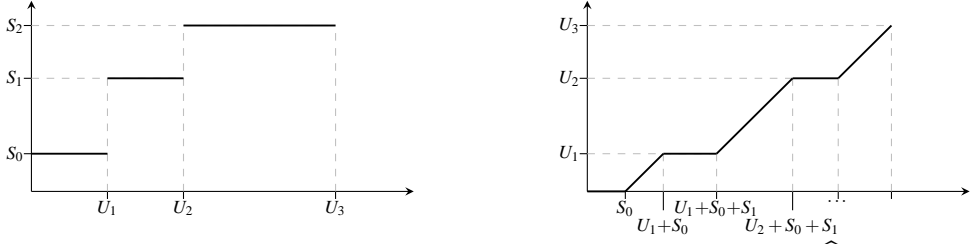

\begin{minipage}{0.49\linewidth}
\centering
\begin{pgfpicture}{0mm}{0mm}{50mm}{25mm}
  \pgfsetxvec{\pgfpoint{1mm}{0mm}}
  \pgfsetyvec{\pgfpoint{0mm}{1mm}}
  \pgfsetendarrow{\pgfarrowsingle}
  \pgfline{\pgfxy(0,0)}{\pgfxy(0,25)}
  \pgfline{\pgfxy(0,0)}{\pgfxy(50,0)}
  \pgfsetendarrow{}
  \pgfline{\pgfxy(10,0)}{\pgfxy(10,-1)}
  \pgfline{\pgfxy(20,0)}{\pgfxy(20,-1)}
  \pgfline{\pgfxy(40,0)}{\pgfxy(40,-1)}
  \pgfputat{\pgfxy(10,-1.2)}{\pgfbox[center,top]{\Tiny $U_1$}}
  \pgfputat{\pgfxy(20,-1.2)}{\pgfbox[center,top]{\Tiny $U_2$}}
  \pgfputat{\pgfxy(40,-1.2)}{\pgfbox[center,top]{\Tiny $U_3$}}
  \pgfline{\pgfxy(0,5)}{\pgfxy(-1,5)}
  \pgfline{\pgfxy(0,15)}{\pgfxy(-1,15)}
  \pgfline{\pgfxy(0,22)}{\pgfxy(-1,22)}
  \pgfputat{\pgfxy(-1,5)}{\pgfbox[right,center]{\Tiny $S_0$}}
  \pgfputat{\pgfxy(-1,15)}{\pgfbox[right,center]{\Tiny $S_1$}}
  \pgfputat{\pgfxy(-1,22)}{\pgfbox[right,center]{\Tiny $S_2$}}
  \pgfsetcolor{lightgray}
  \pgfsetdash{{1mm}}{0pt}
  \pgfline{\pgfxy(0,15)}{\pgfxy(20,15)}
  \pgfline{\pgfxy(0,22)}{\pgfxy(40,22)}
  \pgfline{\pgfxy(10,0)}{\pgfxy(10,15)}
  \pgfline{\pgfxy(20,0)}{\pgfxy(20,22)}
  \pgfline{\pgfxy(40,0)}{\pgfxy(40,22)}
  \pgfsetcolor{black}
  \pgfsetlinewidth{0.25mm}
  \pgfsetdash{}{0pt}
  \pgfline{\pgfxy(0,5)}{\pgfxy(10,5)}
  \pgfline{\pgfxy(10,15)}{\pgfxy(20,15)}
  \pgfline{\pgfxy(20,22)}{\pgfxy(40,22)}
\end{pgfpicture}
\end{minipage}
\begin{minipage}{0.49\linewidth}
\centering
\begin{pgfpicture}{0mm}{0mm}{50mm}{25mm}
  \pgfsetxvec{\pgfpoint{1mm}{0mm}}
  \pgfsetyvec{\pgfpoint{0mm}{1mm}}
  \pgfsetendarrow{\pgfarrowsingle}
  \pgfline{\pgfxy(0,0)}{\pgfxy(0,25)}
  \pgfline{\pgfxy(0,0)}{\pgfxy(50,0)}
  \pgfsetendarrow{}
   \pgfline{\pgfxy(5,0)}{\pgfxy(5,-1)}
   \pgfline{\pgfxy(10,0)}{\pgfxy(10,-2.5)}
   \pgfline{\pgfxy(17,0)}{\pgfxy(17,-1)}
   \pgfline{\pgfxy(27,0)}{\pgfxy(27,-2.5)}
   \pgfline{\pgfxy(33,0)}{\pgfxy(33,-1)}
   \pgfline{\pgfxy(40,0)}{\pgfxy(40,-1)}
   \pgfputat{\pgfxy(5,-1.2)}{\pgfbox[center,top]{\Tiny $S_0$}}
   \pgfputat{\pgfxy(10,-3.6)}{\pgfbox[center,top]{\Tiny $U_1$$+$$S_0$}}
   \pgfputat{\pgfxy(17,-1.2)}{\pgfbox[center,top]{\Tiny $U_1$$+$$S_0$$+$$S_1$}}
    \pgfputat{\pgfxy(27,-3.6)}{\pgfbox[center,top]{\Tiny $U_2+S_0+S_1$}}
    \pgfputat{\pgfxy(33,-1.2)}{\pgfbox[center,top]{\Tiny $\dots$}}
  \pgfline{\pgfxy(0,5)}{\pgfxy(-1,5)}
  \pgfline{\pgfxy(0,15)}{\pgfxy(-1,15)}
  \pgfline{\pgfxy(0,22)}{\pgfxy(-1,22)}
  \pgfputat{\pgfxy(-1,5)}{\pgfbox[right,center]{\Tiny $U_1$}}
  \pgfputat{\pgfxy(-1,15)}{\pgfbox[right,center]{\Tiny $U_2$}}
  \pgfputat{\pgfxy(-1,22)}{\pgfbox[right,center]{\Tiny $U_3$}}
  \pgfsetcolor{lightgray}
  \pgfsetdash{{1mm}}{0pt}
  \pgfline{\pgfxy(0,5)}{\pgfxy(17,5)}
  \pgfline{\pgfxy(0,15)}{\pgfxy(27,15)}
  \pgfline{\pgfxy(0,22)}{\pgfxy(40,22)}
  \pgfline{\pgfxy(10,0)}{\pgfxy(10,5)}
  \pgfline{\pgfxy(17,0)}{\pgfxy(17,5)}
  \pgfline{\pgfxy(27,0)}{\pgfxy(27,15)}
  \pgfline{\pgfxy(33,0)}{\pgfxy(33,15)}
  \pgfline{\pgfxy(40,0)}{\pgfxy(40,22)}
  \pgfsetcolor{black}
  \pgfsetlinewidth{0.25mm}
  \pgfsetdash{}{0pt}
  \pgfline{\pgfxy(0,0)}{\pgfxy(5,0)}
  \pgfline{\pgfxy(5,0)}{\pgfxy(10,5)}
  \pgfline{\pgfxy(10,5)}{\pgfxy(17,5)}
  \pgfline{\pgfxy(17,5)}{\pgfxy(27,15)}
  \pgfline{\pgfxy(27,15)}{\pgfxy(33,15)}
  \pgfline{\pgfxy(33,15)}{\pgfxy(40,22)}
\end{pgfpicture}
\end{minipage}
\vspace{2mm}
\caption{A picture of $D_t$ (left) and of $E_t$ (right) when $\Xz_0=0$.}\label{f:de}
\end{figure}

Finally, we give an explicit construction (which is taken from Example $6.31$
of~\cite{EngSch85}). Denote by $L_0$ the local time of $\Xz$ in $0$ and
consider $(S_n, S_n')$ independent exponential random variables of rate $\lambda$,
independent from $\Xz$. Define $U_0=0$ and
$$
U_{n+1} = U_n + \inf\{t\geq0: L_0(t+U_n) - L_0(U_n) > S_n'\}.
$$
Define finally the time-delay $(D_t)_{t\geq0}$ as
$$
D_t = S_0\uno_{\{\Xz_0=0\}} + \sum_{k=1}^\infty S_k\uno_{[U_k,+\infty)}(t)
$$
and denote by $(E_t)_{t\geq0}$ the inverse of $t+D_t$ (see Figure~\ref{f:de}).
The process $Y_t^x = \Xz_{E_t}$ is a Markov process.

\noindent In few words, $D_t$ jumps (and so $Y_t^x$ stops at zero for an amount of time
corresponding to the size $S$ of the jump) every time the local time $L_0$ accumulates
enough \emph{mass} (in terms of random variables $S'$).
\begin{remark}
A similar example can be considered in dimension $2$ (or more), but the
behaviour of solutions is slightly different, see Example $4.12$ of
Stroock \& Yor~\cite{StrYor80}.
\end{remark}
\subsubsection{Analysis of invariant measures in a dumped version}

Similar conclusions can be drawn for the damped problem,
$$
d\tX_t = -\tX_t\,dt + \sigma_\alpha(\tX_t)\,dW_t.
$$
Indeed, we can use the method of \emph{removal of drift} (see for
example Proposition 5.13 in Karatzas \& Shreve~\cite{KarShr91}),
and reduce the problem to an equation of the same type of~\eqref{e:girsanov},
with a different diffusion coefficient, which anyway has exactly the same
regularity properties as $\sigma_\alpha$. This is possible since
$\tfrac{b(x)}{\sigma_\alpha(x)^2}$ is a bounded function, where
$b$ is the drift function $b(x)=-x$.

It is easy to verify that each Markov solution has a unique invariant
measure. Each of these measures, with the exception of the one corresponding
to the Markov process which spends no time in $0$, has an atom in $0$.
In particular, there are invariant measures that are not mutually
equivalent.
\subsection{A \emph{strong Feller} example by Stroock and Yor}\label{ss:stroockyor}

Following Example 4.5 of Stroock \& Yor~\cite{StrYor80}, consider
the following diffusion operator,
$$
\genSV = \frac12\uno_G(x)\frac{\partial^2}{\partial x^2} + \uno_{\{0\}}(x)\frac{\partial}{\partial x},
$$
where $G=\Erre\backslash\{0\}$, and denote by $\msSV(x)$, for every
$x\in\Erre$, the set of all probability measures on $C([0,\infty);\Erre)$
solutions to the martingale problem associated to $\genSV$.

For every $x\in\Erre$, denote by $W_x$ the law of $x+B_t$, where
$(B_t)_{t\geq0}$ is a one-dimensional standard Brownian motion
(hence $W_x$ is the Wiener measure at $x$) and it is clear that
$W_x\in\msSV(x)$ for all $x$. In particular, $(W_x)_{x\in\Erre}$
is a Markov solution to the problem which is \emph{strong Feller}.

The problem is not well posed and it is possible to see that
there is another strong Feller Markov solution, corresponding
to the reflected Brownian motion. We give a few hints, all
details can be found in Stroock \& Yor~\cite[Example 4.5]{StrYor80}.

First, by Lemma 4.6 of~\cite{StrYor80}, a probability measure
$P\in\msSV(x)$ if and only if
\begin{enumerate}
\item $P[\xi_0=x]=1$,
\item for every $\phi\in C^{1,2}([0,\infty)\times\Erre)$ such that
      $\partial_t\phi(t,0) + \partial_x\phi(t,0)\geq0$,
      $$
      \phi(t,\xi_t) - \int_0^t\uno_G(\xi_s)\bigl[\partial_t\phi(s,\xi_s) + \frac12\partial_x^2\phi(s,\xi_s)\bigr]\,ds
      $$
      is a $P$-submartingale with respect to the natural filtration $(\B_t)_{t\geq0}$.
\end{enumerate}
By Theorems $3.1$ and $5.5$ of Stroock \& Varadhan~\cite{StrVar71}, for every
$x\geq0$ there exists a unique probability measure $Q_x$ such that $Q_x\in\msSV(x)$ 
and
\begin{equation}\label{e:Qunique}
Q_x[\xi_t\geq0\text{ for all }t\geq0] = 1.
\end{equation}
Define $\tau_0$ as the hitting time of $0$. If $x<0$, define
$Q_x$ as the probability measure equal to $W_x$ up to time
$\tau_0$, and then equal to $Q_0$ suitably translated to
time $\tau_0$ afterwards. Property~\eqref{e:Qunique}
ensures that the solution $(Q_x)_{x\in\Erre}$ is Markov.
Moreover, since $Q_x=W_x$ up to time $\tau_0$, it follows
that 
$$
\E^{Q_x}[\phi(\xi_t)]
 = \E^{W_x}[\phi(\xi_t)\uno_{\{\tau_0>t\}}] + \E^{W_x}[\widetilde\semi_{t-\tau_0}\phi(0)\uno_{\{\tau_0\leq t\}}],
$$
for every $\phi\in C_b(\Erre)$, where $(\widetilde\semi_t)_{t\geq0}$ is the transition
semigroup associated to $(Q_x)_{x\in\Erre}$. Hence $(Q_x)_{x\in\Erre}$ is also strong
Feller.
\end{document}